\documentclass[preprint,12pt]{elsarticle}

\usepackage[english]{babel}

\usepackage[a4paper, top=2cm,bottom=2cm,left=3cm,right=3cm,marginparwidth=1.75cm]{geometry}

\usepackage{amsmath}
\usepackage{amssymb}
\usepackage{graphicx}
\usepackage{algorithm}
\usepackage{algpseudocode}
\usepackage{listings}
\usepackage[colorlinks=true, allcolors=blue]{hyperref}
\usepackage{float}
\usepackage{enumerate}
\usepackage{multirow}
\usepackage{caption}
\usepackage{subcaption}
\usepackage{cancel}

\def\FZcomment[#1]{\textcolor{blue}{#1}}
\newcommand{\name}{{JA\_par}}
\begin{document}
\begin{frontmatter}
\title{Linearizing Anhysteretic Magnetization Curves: A Novel Algorithm for Finding Simulation Parameters and Magnetic Moments}
\author[a]{Daniele Carosi} 
\author[b]{Fabiana Zama}
\author[a]{Alessandro Morri}
\author[a]{Lorella Ceschini}

\affiliation[a]{organization={Department of Industrial Engineering, University of Bologna},
            country={Italy}}
\affiliation[b]{organization={Department of Mathematics, University of Bologna},
            country={Italy}}

\begin{abstract}
This paper proposes a new method for determining the simulation parameters of the Jiles-Atherton Model used to simulate the first magnetization curve and hysteresis loop in ferromagnetic materials. The Jiles-Atherton Model is an important tool in engineering applications due to its relatively simple differential formulation. However, determining the simulation parameters for the anhysteretic curve is challenging. Several methods have been proposed, primarily based on mathematical aspects of the anhysteretic and first magnetization curves and hysteresis loops. This paper focuses on finding the magnetic moments of the material, which are used to define the simulation parameters for its anhysteretic curve. The proposed method involves using the susceptibility of the material and a linear approximation of a paramagnet to find the magnetic moments. The simulation parameters can then be found based on the magnetic moments. The method is validated theoretically and experimentally and offers a more physical approach to finding simulation parameters for the anhysteretic curve and a simplified way of determining the magnetic moments of the material.
\end{abstract}

\end{frontmatter}
\begin{flushleft}

\section{Introduction}
Ferromagnetic materials have long presented a challenge in determining their magnetic constitutive laws. Numerous approaches and mathematical models have been developed to address this issue. The most accurate models, according to the literature, are the Brillouin and Langevin Functions for describing reversible magnetic transformations, which produce "anhysteretic curves," and the Presiach and Jiles-Atherton Model for describing irreversible magnetic transformations, which produce the first magnetization curve and hysteresis loop \cite{Fiorillo2006}.

In daily applications, the magnetization of a magnetic material due to an external generated magnetic field does not pass through equilibrium states but through non - equilibrium states, showing the phenomenon of hysteresis. \\
Many models to describe the hysteresis behaviour of a magnetic material have been developed in the years, like Preisach Model, Stoner-Wolfhart Model and so on. One of the most used, especially in engineering applications is the Jiles-Atherton Model (JA) \cite{Fiorillo2006}.\\
This model describes the magnetisation of a ferromagnetic material in function of the external applied field with a first-order Ordinary Differential Equation (ODE) depending on several critical parameters related to the material and experiment conditions.
To define the JA model of a given material, it is necessary to estimate the parameters from magnetisation measurements at different intensities of the applied magnetic field. Such a problem is a well-known, very difficult task, and many approaches can be found in the literature  
 to simulate the first magnetization curve and hysteresis loop. One such method is the genetic algorithm, which uses a penalty fitness function and boundary values \cite{Chwastek2006}. Another method is the "Branch and bound method," which uses boundary conditions of the parameters and is mainly based on mathematical considerations \cite{Chwastek2007}. A third method involves considering the anhysteretic function similar to the first magnetization curve at the maximum applied field \cite{Hauser_2009}. More recently, an improved genetic algorithm has been developed that uses a loss function to evaluate the distance between simulated and experimental hysteresis loops \cite{Liu2020}. Additionally, neural networks have been used, with inputs such as frequency, maximum flux density and flux density, and a parameter indicating whether the magnetic field increases or decreases \cite{Tian2021}.\\
The above methods determine the simulation parameters primarily based on mathematical aspects of the anhysteretic and first magnetization curves and hysteresis loops, using differential or non-linear equations and differential susceptibilities. However, the parameters for simulating the anhysteretic curve are related to the magnetic moment $m$, which has yet to be determined.\\

The present work aims to investigate a robust way to define approximate parameter values using the material physical properties by introducing the linearization of the anhysteretic Magnetization curves.
This method offers several benefits, such as finding the magnetic moments of the material, finding the simulation parameters for the anhysteretic curve more physically, and finding the simulation parameters based solely on the value of initial anhysteretic susceptibility. The results show that it is possible to describe the anhysteretic magnetization curve of a ferromagnetic material  with a paramagnetic function linearly approximated for every value of the external applied field.
This approach could also be used to define the starting guess of parameter estimation procedure,  making it
robust and efficient.

The paper has the following structure: in section \ref{sec:1} we present the problem and the algorithm in paper \cite{Jiles1992}. Tn section \ref{sec2} we present our estimation procedure and finally in section \ref{sec3} we validate the proposed algorithm.
\section{The problem}\label{sec:1}
According to the JA model, the magnetisation $M$ of ferromagnetic materials in function of an externally applied field $H$ is described by the following ODE:
\begin{equation}\label{eq1}
\frac{dM}{dH}=\frac{1}{1+c} \frac{M_{an}(H)-M}{\delta k - \alpha (M_{an}(H)-M)}+\frac{c}{1+c} \frac{dM_{an}(H)}{dH}
\end{equation}
where $M_{an}(H)$ is the anhysteretic magnetisation function:
\begin{equation} \label{eq2}
M_{an}(H)=M_s \left(\coth \left(\frac{H+\alpha M}{a_J} \right)-\frac{a_J}{H+\alpha M} \right) 
\end{equation}
and $ a_J, c, \delta, k, \alpha$ are   the model parameters \cite{Jiles1986}, defined in table \ref{tab:JA}.

\begin{table}[h!]
\begin{tabular}{c|l}
\hline
$ a_J=k_B T/(\mu_0 m) $ & related to the shape of the anhysteretic curve\\
$k_B$ & the Boltzmann constant\\
$T$ & the temperature of the material in $K$\\
$\mu_0$ & the magnetic permeability of free space\\
$m$ & the magnetic moment of a pseudo-domain\\
$\alpha$ & related to the interdomain coupling\\
$M_s$ & the saturation magnetisation\\
$k$ & related to the coercive field and the pinning sites\\
$c$ & related to the reversible processes of magnetisation\\
$\delta=sign(dH/dt)$ & related to the derivative of external applied magnetic field \\
\hline
\end{tabular}
\caption{JA model parameters.}
\label{tab:JA}
\end{table}


One of the most difficult tasks of such a modelling problem is the determination of the model  parameters from measurements of the anhysteretic magnetization $M_{an}$ at different intensities of the external field $H$. 
In addition to the inherent difficulty of solving an ill-posed problem, JA model also presents extreme challenges in defining starting  guesses and extreme sensitivity to their value. 

To address this difficulty, we propose to improve the original 
the approach of  \cite{Jiles1992} based on the exploitation of physical relationships between different quantities that can be obtained from the measurements.
Such quantities are reported in table \ref{tab:mJA}
\begin{table}[h!]
\begin{tabular}{c|l}
\hline
$\chi'_{in} $ &  Initial differential susceptibility of the first magnetisation curve \\
$\chi'_{an}$ &  Initial differential susceptibility of the anhysteretic magnetisation curve \\
$\chi'_{max}$ & the differential susceptibility at coercive field \\
$\chi'_r$ & Differential susceptibility at remanence point \\
$\chi'_m$ & Differential susceptibility at hysteresis loop tip \\
$H_c$ & Value of coercive field \\
$M_r$ & Value of magnetisation at remanence point \\
$M_m$ & Value of magnetisation at loop tip \\
$H_m$ & Value of external applied field corresponding to $M_m$\\
\hline
\end{tabular}
\caption{Physical quantities that can be extracted from measured data.}
\label{tab:mJA}
\end{table}
The estimation procedure proposed in \cite{Jiles1992} exploits the quantities reported in table \ref{tab:mJA} and reduces the dependence of all model parameters to a single parameter $\alpha$ which is heuristically set. Such procedure is outlined in algorithm \ref{alg:estimation}.
\begin{algorithm}[h]
\caption{Estimation Algorithm Jiles \cite{Jiles1992}}\label{alg:estimation}
\begin{algorithmic}[1]
\State Calculate the value of $c$ with equation: $c=\frac{\chi'_{in}}{\chi'_{an}}$
\Repeat
\State Set a seed value of $\alpha$
\State Calculate a first estimation value of $a_J$ with equation: $a_J = \frac{M_s}{3} \left( \frac{1}{\chi' _{an}} + \alpha \right)$
\State Compute $k$ as:
$$k=\frac{M_{an}(H_c)}{1-c} \left (\alpha + \frac{1}{ \left (\frac{1}{1-c} \right) \chi'_{max} - \left (\frac{c}{1-c} \right) \frac{dM_{an}(H_c)}{dH} } \right) $$
\State Solve for $\alpha$  the following non-linear equation, using the current $\alpha$ estimate as  initial guess: 
$$M_r = M_{an}(M_r) + \frac{k}{ \left (\frac{\alpha}{1-c} \right) + \frac{1}{ \chi'_r - c \frac{dM_{an}(M_r)}{dH} } }$$
\State Update  $a_J$  solving numerically the following non-linear equation, using as initial guess the current value of $a_J$: $$M_m = M_{an}(H_m) - \frac{ (1-c) k \chi'_m}{\alpha \chi'_m +1}$$
\State Solve \ref{eq1} with the estimated parameters in the measured points
\Until {\tt fit\_condition}
\end{algorithmic}
\end{algorithm}
The ${\tt fit\_condition}$ is usually represented by a test on the least squares distance or Mean Square Error  between the  simulated hysteresis loop and the experimental data. A known weakness of such an approach is the lack of guarantee that the 
exit condition can be fulfilled; therefore, a restart with different seeds might be  required, depending on the measured data. 
One strength is the reduced computational cost consisting of the solution of two nonlinear equations for each iteration.

In the next section we introduce a method to find the magnetic moment $m$ based on the linearisation of the anhysteretic magnetisation curve and its susceptibilities. This provides a simple way to find values of parameters $a_J$ and $\alpha$.

\section{ Magnetic moments and simulation parameters of anhysteretic curve} \label{sec2}
Let's start considering an anhysteretic theoretic curve of a ferromagnetic material generated by \eqref{eq2} with given simulation parameters $\alpha$ and $a_J$.\\
Since for every value of external applied field $H_a$ the curve has only one value of magnetisation $M$, the following injective function can describe its behaviour: 
\begin{equation} \label{eq3}
M_{an}(H)=M_s \left(\coth \left(\frac{H}{a} \right)-\frac{a}{H} \right) 
\end{equation}
    with a simulation parameter $a \neq a_J$.  
    Since there is no  interaction between the magnetic moments in paramagnetic materials, the parameter $\alpha$ is absent. 
Injective functions such as \eqref{eq3} usually describe the magnetic behaviour of a paramagnetic material \cite{Bertotti1998}. 
Therefore, analogously to \eqref{eq2}, the shape parameter 
 $a$ is:
\begin{equation}\label{eq:a}
a=\frac{k_B T}{\mu_0 m_1},
\end{equation} where $m_1$ is the magnetic moment of an equivalent paramagnetic curve that can describe the ferromagnetic one. \\
A common experimental procedure to obtain the anhysteretic curve of a magnetic material consists of  superimposing a steady external magnetic field on another magnetic field that varies between a  minimum and maximum value. 
The varying magnetic field is responsible for creating the hysteresis loop. Gradually, the range of the varying magnetic field is reduced until it aligns with the value of the constant external magnetic field. Through this procedure, consistent magnetization values of the material can be obtained that are free of hysteresis \cite{Tumanski2016}.\\
Considering then that the magnetization curve of ferromagnetic materials is obtained with relatively small and constant values of the external applied field  $H_a$, it is possible to obtain a  simplified approximation of $M_{an}(H)$ considering  the Taylor expansion  of equation \eqref{eq3} : 
$$M_{an}(H) \approx M_s \left(\coth \left(\frac{H_a}{a} \right)-\frac{a}{H_a} \right) + \mathcal{O}(H-H_a)$$
and taking only the first term. In the assumption that $H \approx H_a$, we have:
$$
M_{an}(H) \approx M_s \left(\coth \left(\frac{H_a}{a} \right)-\frac{a}{H_a} \right) 
$$
then considering the series expansion of $\coth(x)$ for $x \neq 0$, 
$$\coth(x) = \frac{1}{x} + \frac{x}{3} - \frac{x^3}{45} + \frac{2x^5}{945} - \cdots
$$
and taking the first two terms,  we obtain the linearized approximation $M_{an}^a$ :
$$M_{an}^a = M_s \left ( \cancel{\frac{1}{\tfrac{H_a}{a}}}+ \frac{\tfrac{H_a}{a}}{3}- \cancel{\frac{a}{H_a}} \right )  $$
Substituting $a$ from  \eqref{eq:a}:
 \begin{equation} \label{eq5}
M_{an}^a = \frac{M_s}{3} \frac{\mu _0 m_1}{k_B T} H_a
\end{equation}
Using $M_{an}^a$ to define the anhysteretic susceptibility of the ferromagnetic material $\chi_{an}^a$ \cite{Tumanski2016},\cite{Bertotti1998} i.e.:
\begin{equation} \label{eq6}
\chi_{an}^a = \frac{M_{an}^a}{H_a}
\end{equation}
and substituting into   \eqref{eq5}
it is possible to compute the magnetic moment $m_1$ of the equivalent paramagnetic material for every value of the external applied field   related to the anhysteretic susceptibility of the ferromagnetic material:
\begin{equation} \label{eq30}
m_1 = \frac{3 k_B T \chi_{an}^a}{\mu_0 M_s}.
\end{equation}
Going back to equation \eqref{eq2} that describes the anhysteretic behaviour of a ferromagnetic material we have:
\begin{equation} \label{eq31}
M_{an}= \frac{M_s}{3 a_J}(H_a + \alpha M_{an})=\frac{M_s}{3}  \frac{\mu_0 m}{k_B T} (H_a + \alpha M_{an}) 
\end{equation}
To find the value $m$ for the magnetic moment of the ferromagnetic material in \eqref{eq31}, we substitute $M_{an}$ with $M_{an}^a$:
$$
\begin{aligned}
& M_{an}^a & = & \frac{M_s}{3}  \frac{\mu_0 m}{k_B T} (H_a + \alpha M_{an}^a),  \\
& \chi_{an}^a & = & \frac{M_s}{3}  \frac{\mu_0 m}{k_B T} ( 1 + \alpha \chi_{an}^a), \\
\end{aligned}
$$
and find the following  expression for  $m$:
\begin{equation}\label{eq:32_3}
 m  =  \frac{3 k_B T}{M_s}\frac{\chi_{an}^a}{( 1 + \alpha \chi_{an}^a)} 
\end{equation}
However, since $\alpha $ is unknown this relation cannot be applied to compute $m$. We observe 
that $m$ multiplies the external and the molecular field, so setting $\alpha=0$ 
we obtain that $m$ can be considered as  the magnetic moment of an equivalent paramagnetic curve of the ferromagnetic material, 
with  susceptibility  $\chi_{param}$. 
Hence, using  \eqref{eq30}, we obtain an alternative characterization of $m$ depending on the  unknown susceptibility  $\chi_{param}$ :
\begin{equation} \label{eq33}
m = \frac{3 k_B T \chi_{param}}{\mu_0 M_s}.
\end{equation}
%
An estimate of $\chi_{param}$ can be obtained by  substituting \eqref{eq33} in  equation \eqref{eq3} 
\begin{equation} \label{eq34}
M_{an}= M_s \left (\coth \left (\frac{3 \chi_{param} H_a}{M_s} \right) - \frac{M_s}{3 \chi_{param} H_a} \right) 
\end{equation}
and solving the nonlinear equation for properly chosen  values $M_{an}$ and $H_a$.
The idea is to choose very large values of the  applied external field $H_a^1$ since the molecular field in the paramagnetic case does not act. Still, the saturation of the magnetization is almost reached.
However, since  measuring the anhysteretic magnetisation curve for a very high external applied magnetic field value is impossible, we  use its equivalent paramagnetic curve with magnetic moment $m=m_1$.\\
We can compute the corresponding magnetization value $M_{an_1}(H_a^1)$ as: 
\begin{equation} \label{eq37}
M_{an_1}(H_a^1) = M_s \left (\coth \left(\frac{\mu_0 m_1}{k_B T} H_a^1 \right ) - \frac{k_B T}{\mu_0 m_1 H_a^1} \right) 
\end{equation}
Since in general $M_{an_1} > M_{an}$, we solve equation  \eqref{eq34} with the  magnetization value $M_{an}$  estimated as follows:  $$M_{an} \approx \eta^* M_{an_1}, \quad 0.9 \leq \eta^* <1.$$
Finally, we can simplify computation by considering the anhysteretic susceptibility of the equivalent paramagnetic curve
$$\chi_{an_1}= \frac{M_{an_1}}{H_a^1}$$

and substituting it into \eqref{eq37} we obtain :
\begin{equation} \label{eq38}
\eta^*  \cdot \chi_{an_1} - \frac{M_s}{H_a^1} \left (\coth \left (\frac{3 \chi_{param} H_a^1}{M_s} \right) - \frac{M_s}{3 \chi_{param} H_a^1} \right) = 0.
\end{equation}
After computing $\chi_{param}$ through the  numerical solution of \eqref{eq38}, 
we obtain the magnetic moment $m$ through \eqref{eq33} and the parameter 
$\displaystyle a_J=\frac{k_B T}{\mu_0 m}$, then
we can compute  $\alpha$ by relations  \eqref{eq:32_3} and \eqref{eq33}, i.e.:
\begin{equation} \label{eq36}
\alpha  = \frac{1}{\chi_{param}} - \frac{1}{\chi_{an}^a}.
\end{equation}
Once the parameters have been computed, we compute the approximate anhysteretic magnetization value $\hat M_{an}$  correspondent to the observed applied external field $H_a$, solving the following nonlinear equations for each observed field $H_a$ :
\begin{equation}\label{eq:est}
\hat M_{an} - M_s \left(\coth \left(\frac{H_a +\alpha \hat M_{an}}{a_J} \right)-\frac{a_J}{H_a +\alpha \hat M_{an}} \right) = 0.
\end{equation}
Since $\eta^*$ is not known, the idea is to evaluate \eqref{eq38} in  a sequence  $\{\eta_k\}_{k>0} \in [0.9, 1)$ and define
$$\eta^* = \arg\min_{\eta} \|\mathbf{r}^{(\eta)} \|_2, \quad \mathbf{r}^{(\eta)}= \mu_0 \| \hat M_{an} - M_{an} \|_2$$
where $\mathbf{r}^{(\eta)}$ has  components $r_i^{(\eta)}$ given by  the difference between the  magnetic anhysteretic data  $(\mu_0  M_{an_i})$ and its approximation $(\mu_0 \hat M_{an_i})$ for each corresponding value of external applied field.

These steps are is summarized  in Algorithm \name.

\begin{algorithm}[h!]
\caption{Algorithm \name{}. INPUT: $(H_a,  M_{an})$, $k_B$, $T$, $M_s$}\label{alg:num_algo}
\begin{algorithmic}[1]
\State $\chi_{an}^a=\frac{ M_{an}^a(1)}{H_a(1)}$; with $H_a(1) \neq 0$, s.t. $M_{an}^a(1)$ first experimental positive value;
\State $\displaystyle m_1=\frac{3 k_B T \chi_{an}^a}{\mu_0 M_s}$;
\State $\displaystyle a_1=\frac{k_B T}{\mu_0 m_1}$;
\State ${H_a}^1 =10^6 $: 
$\displaystyle M_{an_1}= \left [M_s \left (\coth \left(\frac{H_a^1}{a_1} \right ) - \frac{a_1}{ H_a^1} \right) \right ]$
\State $\displaystyle \chi_{an_1}=\frac{M_{an_1}}{{H_a^1}}$;
\State $\varepsilon=10^{-5}$;\Comment{step increment}
\State $\eta_0=0.9$;
\State $k=0$
\State  {\tt Loop = true}
\While{\tt Loop}
\State Compute  $\chi_{param}$ by solving the nonlinear equation : $$\displaystyle \eta_k \cdot \chi_{an_1} - \frac{M_s}{H_a^1} \left (\coth \left (\frac{3 \chi_{param} H_a^1}{M_s} \right) - \frac{M_s}{3 \chi_{param} H_a^1} \right) = 0; $$ 
\State $m_2=\frac{3 k_B T \chi_{param}}{\mu_0 M_s}$;
\State $a_J=\frac{k_B T}{\mu_0 m_2}$;
\State $\alpha= \frac{1}{\chi_{param}} - \frac{1}{\chi_{an}^a}$;
\State Solve the nonlinear system for $M$: 
$$\displaystyle \hat M_{an} - M_s \left(\coth \left(\frac{H_a +\alpha \hat M_{an}}{a_J} \right)-\frac{a_J}{H_a +\alpha \hat M_{an}} \right) = 0$$
\State $r=\mu_0( \hat M_{an}-  M_{an})$;\Comment{Residual vector}
\State $k=k+1$; $Nr(j)=norm(r)$; \Comment{Norm of residual vector}
\State ${\tt Loop}= Nr(j)<Nr(j-1)$
\If{${\tt Loop}$}
\State $\eta_{k+1}=\eta_k+\varepsilon$
\EndIf
\EndWhile
\end{algorithmic}
\end{algorithm}
From extended experimental tests with different materials it is verified that $\eta \in [0.9, 1)$ is sufficient obtain close ferromagnetic and equivalent paramagnetic curves, for very high value of external applied field.



%
\section{Methodology validation and testing}\label{sec3}
This section validates the proposed method by investigating its robustness in presence of perturbations. Then the method is tested against data taken from literature and real measurements.
\subsection{Methodology  validation}
Firstly, we demonstrate that the anhysteretic curve of a ferromagnetic material can be approximated using the linear approximation of a paramagnet for any given external applied field. By utilizing the anhysteretic susceptibilities of the ferromagnetic material's anhysteretic curve, which are described by equation \eqref{eq6}  and substituted into equation \eqref{eq30},  it becomes possible to approximate the curve for any external applied field value using equation \eqref{eq5} of an equivalent paramagnetic curve. \\
To accomplish this, we generate a theoretical anhysteretic curve by employing real parameters from a ferromagnetic material and synthetic simulation parameters. For example, we consider a carbon steel at room temperature, with the material parameters listed in table \ref{tab:real_par}. We then choose a set of simulation parameters for the JA model, as provided in table 3, to generate an anhysteretic curve, which is depicted in blue in figures \ref{fig:IMG4} and \ref{fig:IMG5_6}.


%
\begin{table}[h!]
\centering
\begin{tabular}{ccc|ccc}
\hline
\multicolumn{3}{c}{Ferromagnetic material parameters} & \multicolumn{2}{c}{JA parameters} \\
\hline
$M_s$ & $T$ & $ T_c$ & $a_J$ & $\alpha$ \\
\hline
\multirow{6}{*}{$1.6 \cdot 10^6 \left [\frac{A}{m} \right]$ }& \multirow{6}{*}{$303.5 [K]$} & \multirow{6}{*}{$1023.5 [K]$ }& $972$ & $1.4 \cdot 10^{-3}$\\
& & & $972$ & $1.0 \cdot 10^{-3}$\\
& & & $972$ & $1.8 \cdot 10^{-3}$\\
& & & $800$ & $1.4 \cdot 10^{-3}$\\
& & & $1000$ & $1.4 \cdot 10^{-3}$\\
& & & $1200$ & $1.4 \cdot 10^{-3}$\\
\hline
\end{tabular}
\caption{Parameters of an electrical steel at room temperature and JA simulation parameters.}
\label{tab:real_par}
\end{table}

%

For each value of the external applied field and magnetization along the anhysteretic curve of the ferromagnetic material, we compute the susceptibility \eqref{eq6} and the values of magnetic moment $m_1$  by applying equation \eqref{eq30}. By substituting these values into equation \eqref{eq5}, we can evaluate the anhysteretic magnetization  for every value of the external applied field.


\begin{figure}[h!]
\centering
\includegraphics[width=0.6\textwidth]{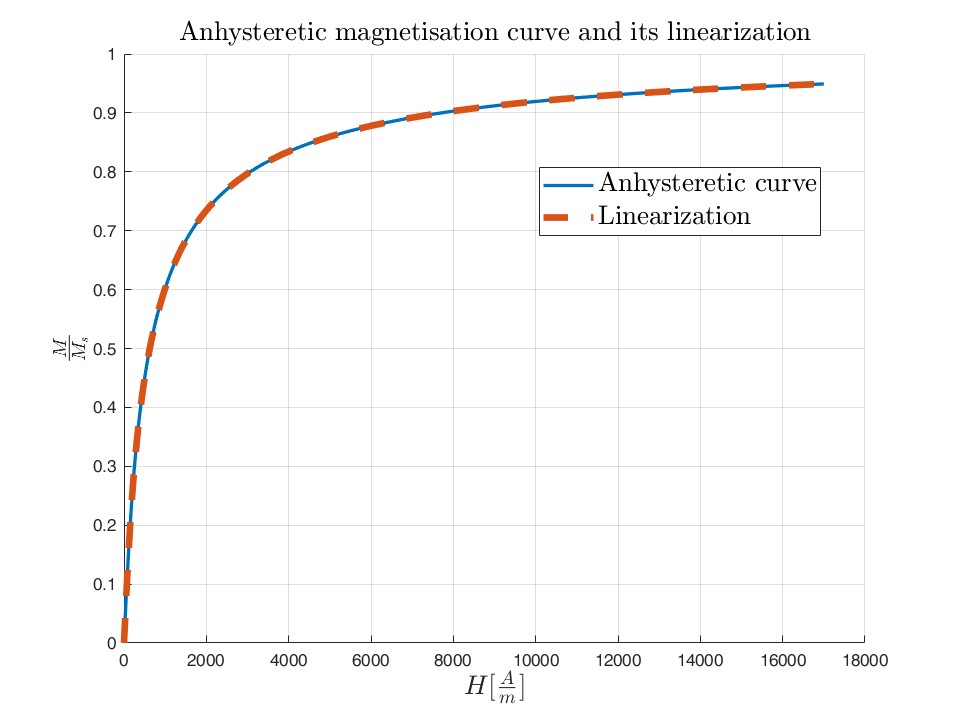}
\caption{\label{fig:IMG1} Theoretical anhysteretic magnetisation curve and its linearization with paramagnetic function.}
\end{figure}
In figure \ref{fig:IMG1} we can appreciate the perfect agreement  between 
the anhysteretic magnetisation curve of a ferromagnetic material and 
that of  an equivalent paramagnetic curve for every value of external applied field.
Such quality is preserved even for changes in the JA parameters $a_J$ and $\alpha$ as reported in the examples represented in figures 
\ref{fig:IMG2_3} where the parameters are modified according to table\ref{tab:real_par} in rows 2-5.
%
%

\begin{figure}[H]
\centering
\includegraphics[height=5cm,width=6cm]{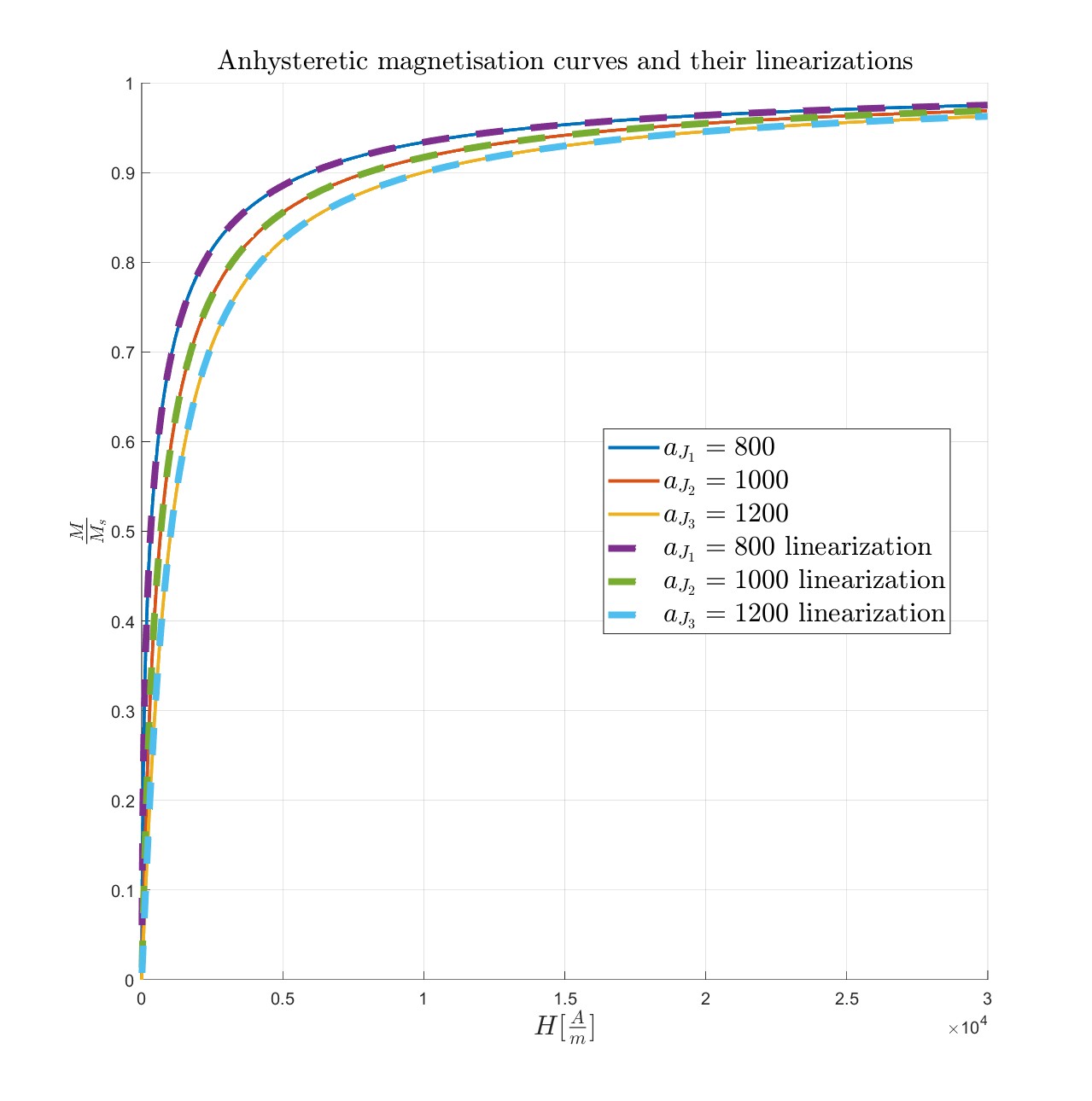}
\includegraphics[height=5cm,width=6cm]{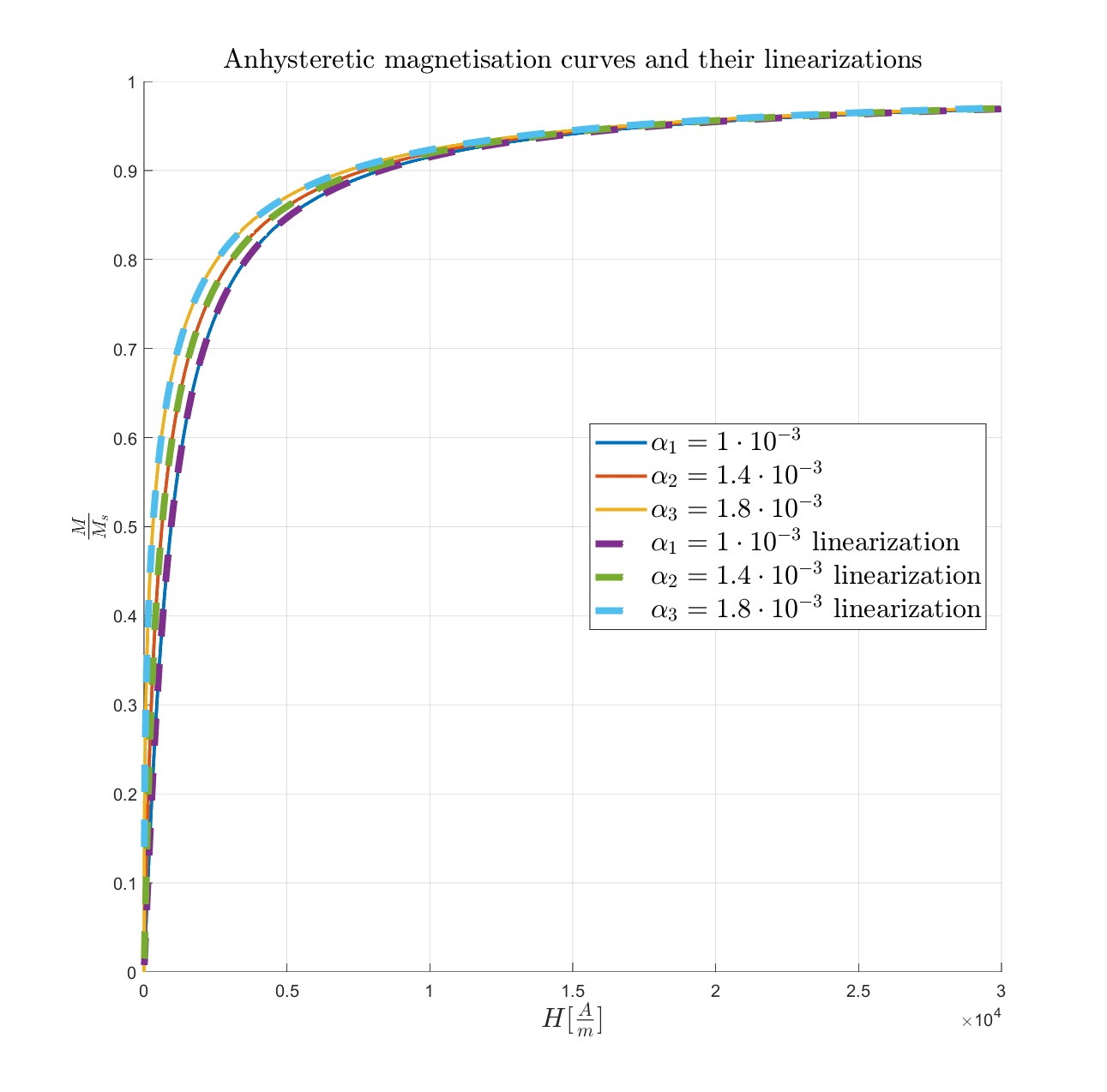}\\
(a) \hspace{3cm} (b)
\caption{\label{fig:IMG2_3} Theoretical anhysteretic magnetisation curves and linearized approximations.
(a) variation of  $a_J$ (b) variation of $\alpha$}
\end{figure}
%
%




%
%
%
By defining the curve obtained using equation \eqref{eq3}, with $a=\frac{\mu_0 m_1}{k_B T}$ where $m_1$ is calculated using the initial anhysteretic susceptibility (considering only the initial value of $H_a>0$), as {\em Paramagnet equivalent 1}, we can observe in Figure \ref{fig:IMG4} that it tends to overestimate the anhysteretic magnetization curve, particularly for small values of $H_a$.

On the other hand, if we consider the curve obtained by setting $\alpha = 0$ in equation \eqref{eq3}, and using the value of $a$ provided in the first row of table \eqref{tab:real_par}, we obtain an underestimating curve, referred to as {\em Paramagnet equivalent 2}. This curve demonstrates better agreement with the anhysteretic magnetization curve, particularly for large values of the applied field $H_a$.



\begin{figure}[h]
\centering
\includegraphics[width=0.8\textwidth]{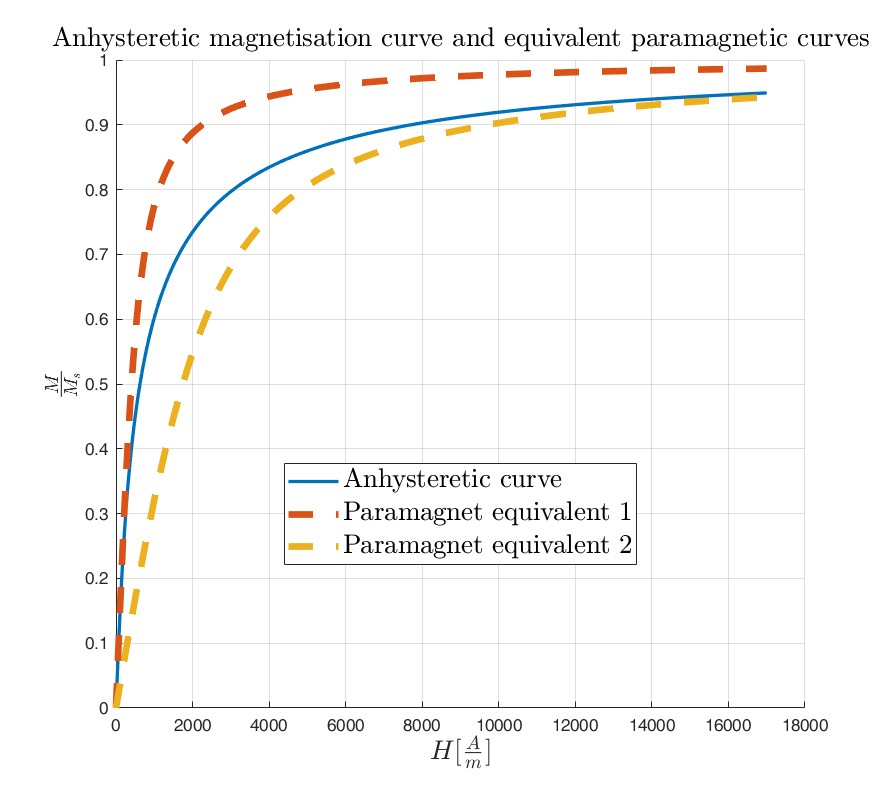}
\caption{\label{fig:IMG4} Theoretical anhysteretic magnetisation curve and its paramagnetic curves}
\end{figure}



Finally, we verify that by evaluating $\chi_{param}$ for extremely high values of external applied field and magnetization, through the solution of equation \eqref{eq:est}, and calculating the initial anhysteretic susceptibility of the ferromagnetic material $\chi_{an}$, we can determine the values of the parameters $a_J$ and $\alpha$ that result in a reliable approximation of the experimental anhysteretic curve of the ferromagnetic material described by equation \eqref{eq2}.
We  validate the evaluation of $\chi_{param}$ and  $\chi_{an}$ by reporting in figures \ref{fig:IMG5_6} the computed magnetization curves varying  $a_J$ and $\alpha$
as in table \ref{tab:real_par}. 

\begin{figure}[h!]
\centering
\includegraphics[height=5cm,width=6cm]{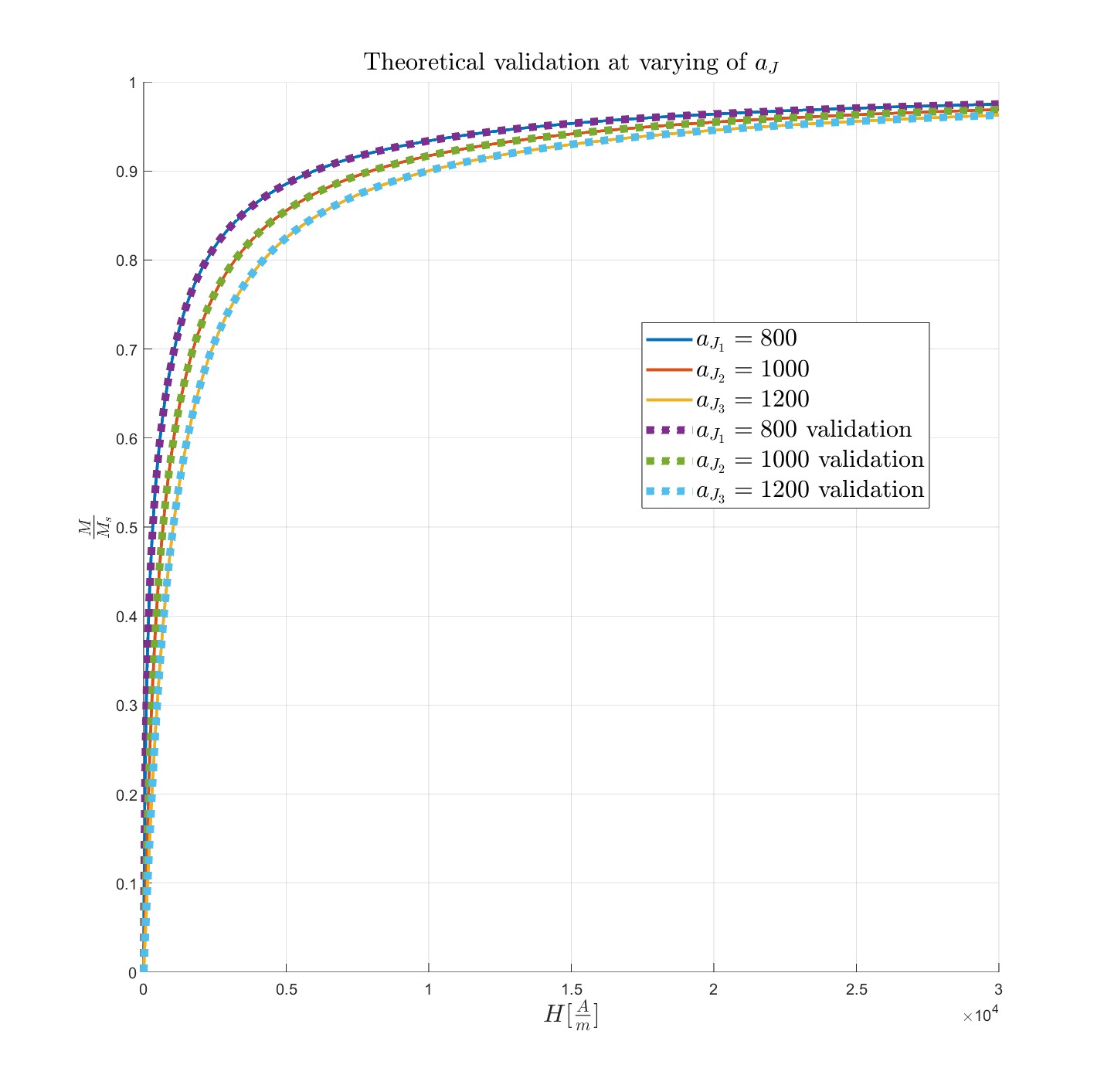}
\includegraphics[height=5cm,width=6cm]{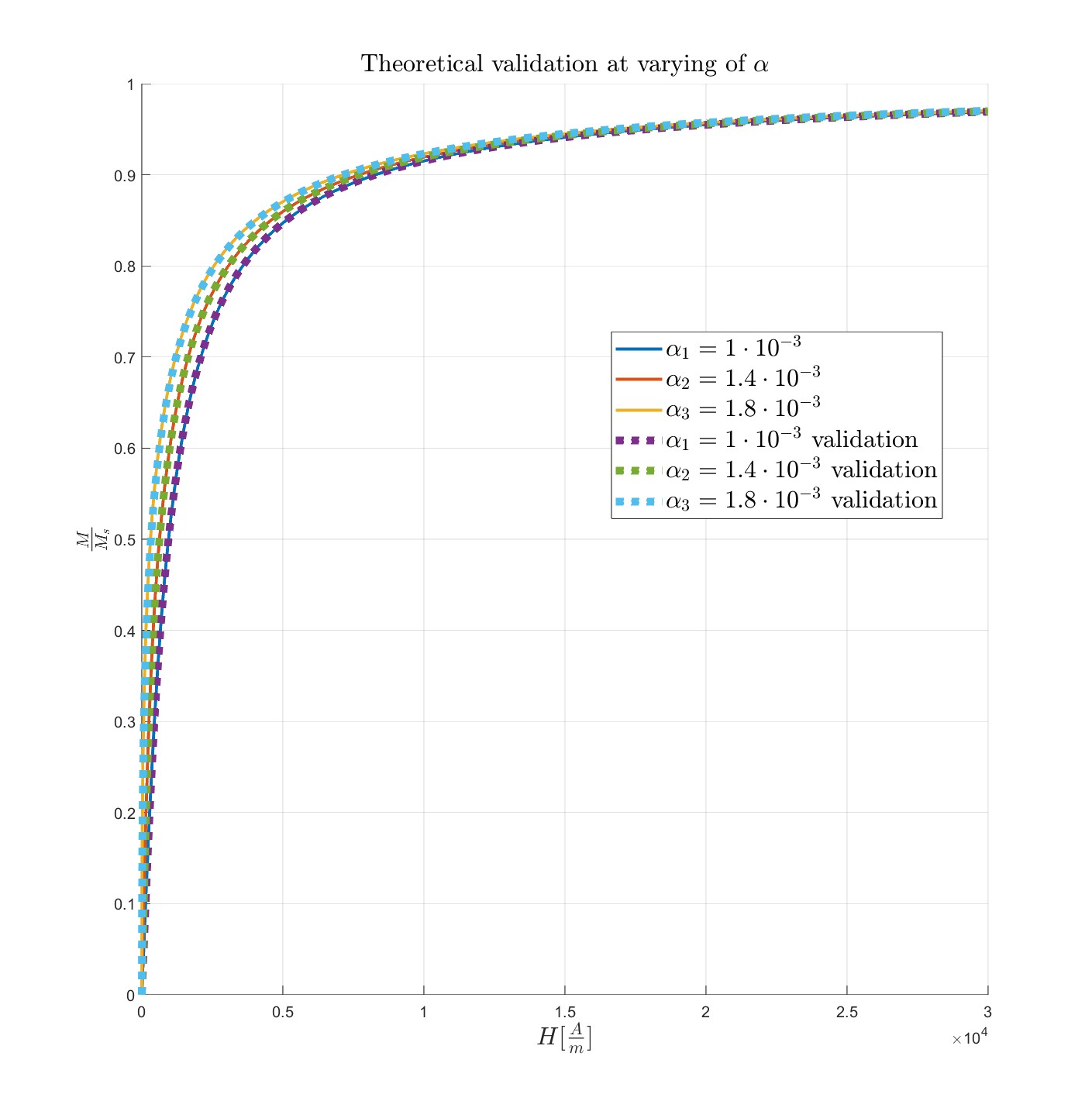}\\
(a) \hspace{3cm} (b)
\caption{\label{fig:IMG5_6} Anhysteretic magnetisation curves  obtaind by $\chi_{param}$, $\chi_{an}$ estimates  and  simulation curves  varying  parameters $\alpha$ and $a_J$. (a) variation of $a_J$  (b) variation of $\alpha$.}
\end{figure}
Again we can observe the perfect agreement between the theoretical and simulated anhysteretic curve with the simulation parameters' variation.

%

\subsection{Algorithm \ref{alg:num_algo} testing}

In this paragraph, we evaluate the performance of the \name{} algorithm using data from papers \cite{Jiles_1984} and \cite{Jiles1986}, which were extracted using the web tool for data extraction called WebPlotDigitilizer \cite{Rohatgi2022}.

Figures \ref{fig:IMG011} depict the residual behaviour within the interval $[0.0, 1]$, thereby confirming that the minimum value can be found in the given interval.
Additionally, these figures provide the optimal value $\eta^*$ computed by the \name{} algorithm.
\begin{figure}[h!]
\centering
\includegraphics[width=6cm,height=5cm]{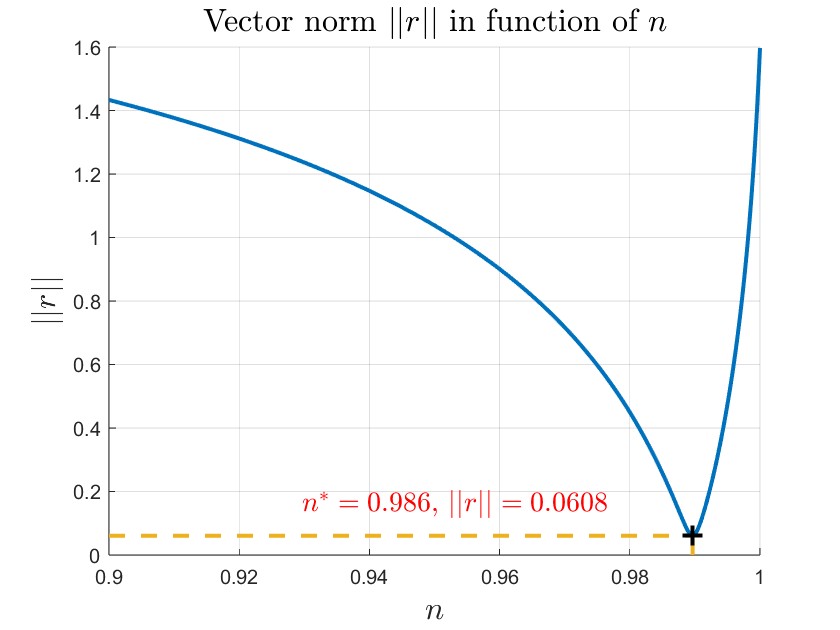}
\includegraphics[width=6cm,height=5cm]{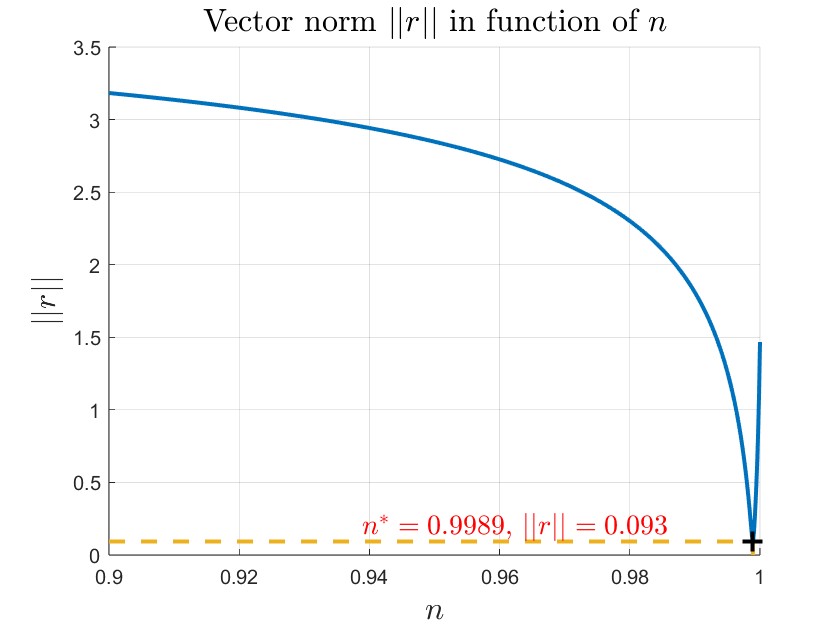} \\
 (a) \hspace{5cm} (b)
\caption{\label{fig:IMG011} Behaviour of the residual norm $|| \textbf{r}||$ for  $n \in [0.9, 1)$. (a) Data in \cite{Jiles_1984} (b) Data in \cite{Jiles1986}.}
\end{figure}
The computed residual and parameters are presented in table \ref{tab:fit1}. 
\begin{table}[h!]
\centering
\begin{tabular}{ c c c c c c c }\hline
 Data & $ \eta^*$ & $|| \textbf{r}||$ & $\chi_{param}$ & $m_2$ & $a_{J}$ & $\alpha$ \\
 \hline
\cite{Jiles_1984} & $0.9896$ & $0.0608$ & $45.7954$ & $2.8738 \cdot 10^{-19}$ & $1.1584 \cdot 10^{4}$ & $0.0195$ \\
\cite{Jiles1986} & $0.9989$ & $0.093$ & $406.5476$ & $2.5512 \cdot 10^{-18}$ & $1.3049 \cdot 10^{3}$ & $0.0021$ \\
\hline
\end{tabular}
\caption{Parameters and residual obtained by algorithm \ref{alg:num_algo}.}
\label{tab:fit1}
\end{table}
The parameters $a_J$ and $\alpha$ corresponding to $\eta^*$ provide the best fit for the anhysteretic data (Figures \ref{fig:IMG6_7}), making them the most representative of the magnetic material.\\
\begin{figure}[h!]
\centering
\includegraphics[height=5cm,width=6cm]{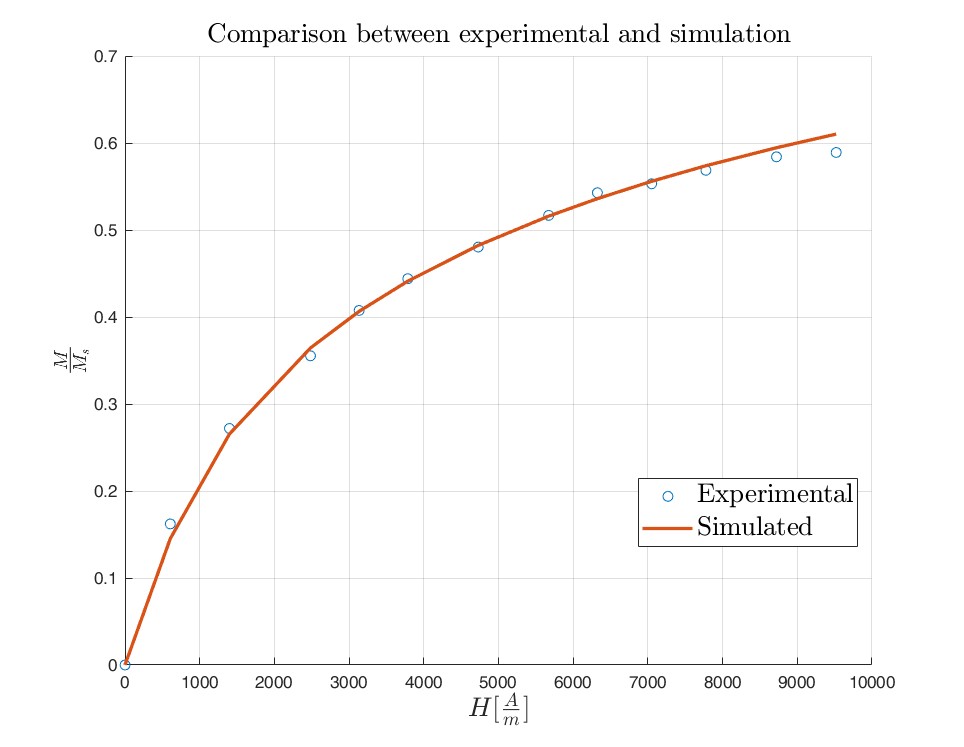}
\includegraphics[height=5cm,width=6cm]{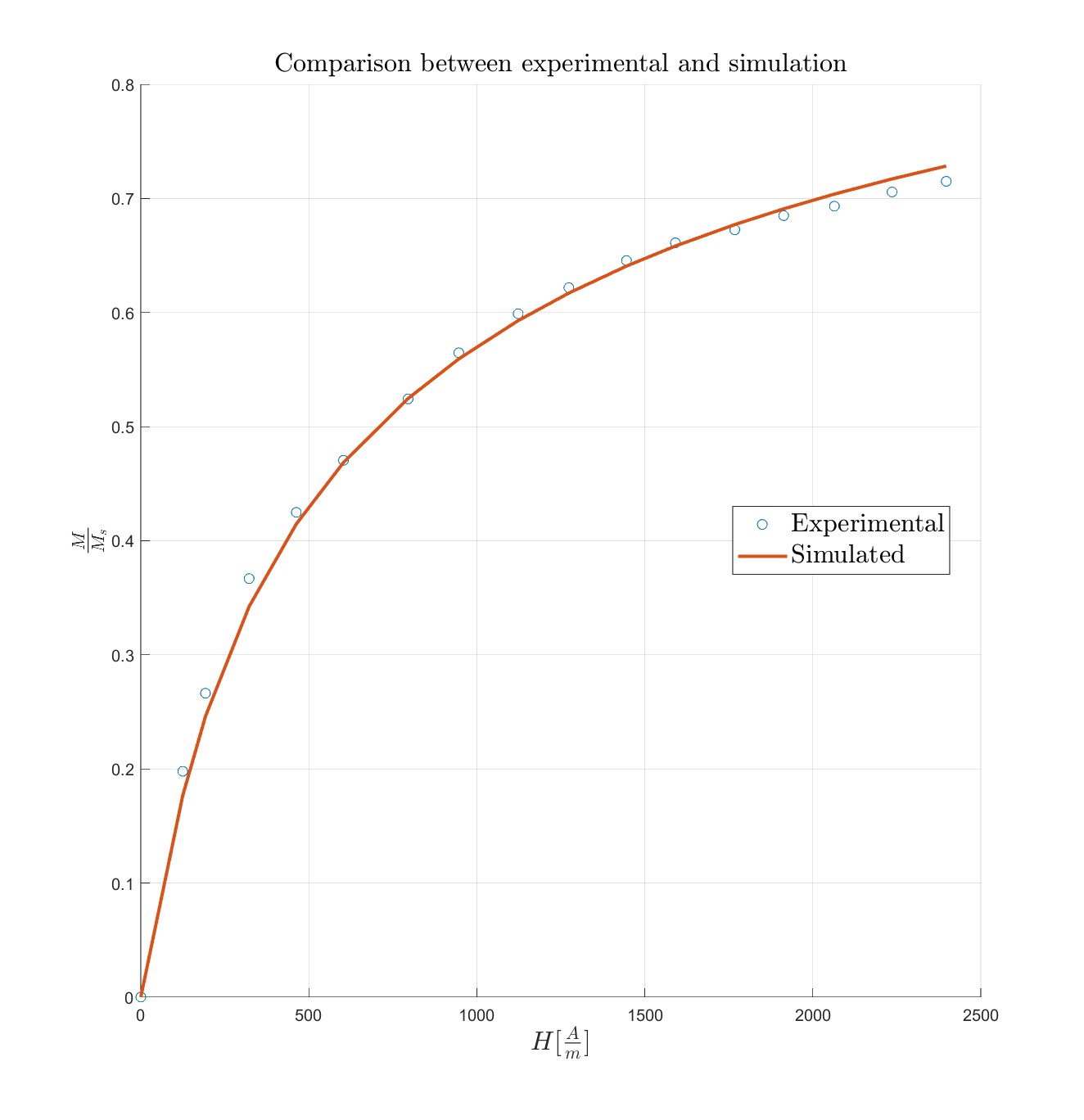}\\
(a) \hspace{3cm} (b)
\caption{Experimental anhysteretic magnetisation curves (blue circles) and  \name{} simulations (red line). (a) data  \cite{Jiles_1984} (b) data  \cite{Jiles1986}}
\label{fig:IMG6_7} 
\end{figure}
From a computational efficiency perspective, we observe that the algorithm requires solving nonlinear equations in steps 11 and 15 of the \name{} algorithm. For this purpose, the function {\tt fzero} is applied, using zero as the starting guess.

\subsection{Experimental hysteresis validation}

It is necessary to verify whether the parameters obtained through simulating the anhysteretic curve and solving equation \eqref{eq38} describe the hysteresis curve accurately. To this purpose, we set the  simulation parameters $c$ and $k$ in \eqref{eq1} as follows:
\begin{itemize}
    \item $c=\frac{\chi'_{in}}{\chi'_{an}}$;
    \item $k = H_c $;
\end{itemize}
where $\chi'_{in}$, $\chi'_{an}$ are defined in table \ref{tab:mJA}.
The results are checked on the curves obtained from Jiles' paper \cite{Jiles1992}  (figures \ref{fig:IMG10_11})  and real measurements (figure \ref{fig:IMG12}).
\begin{figure}[h!]
\centering
\includegraphics[width=7cm,height=7cm]{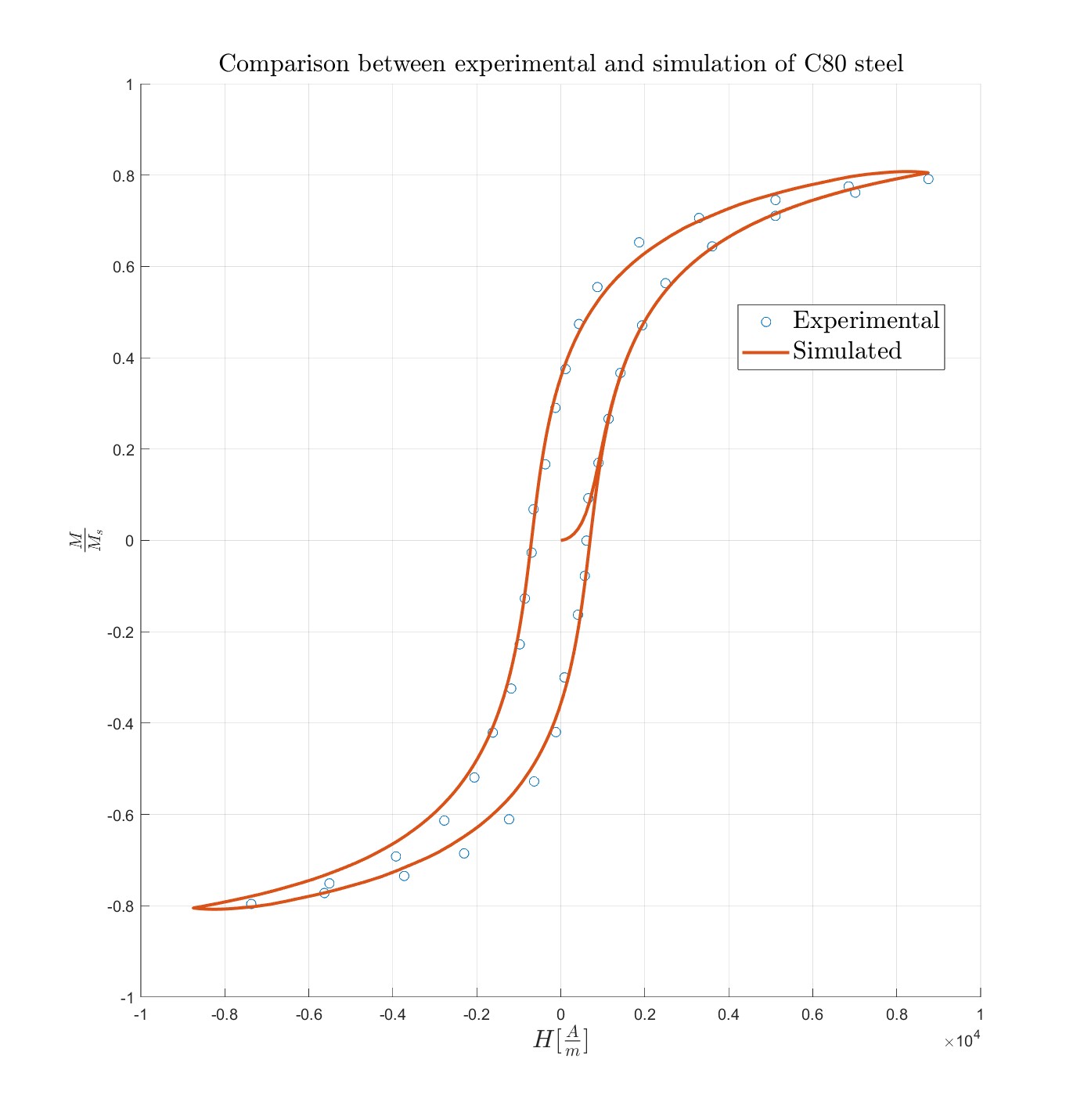}
\includegraphics[width=7cm,height=7cm]{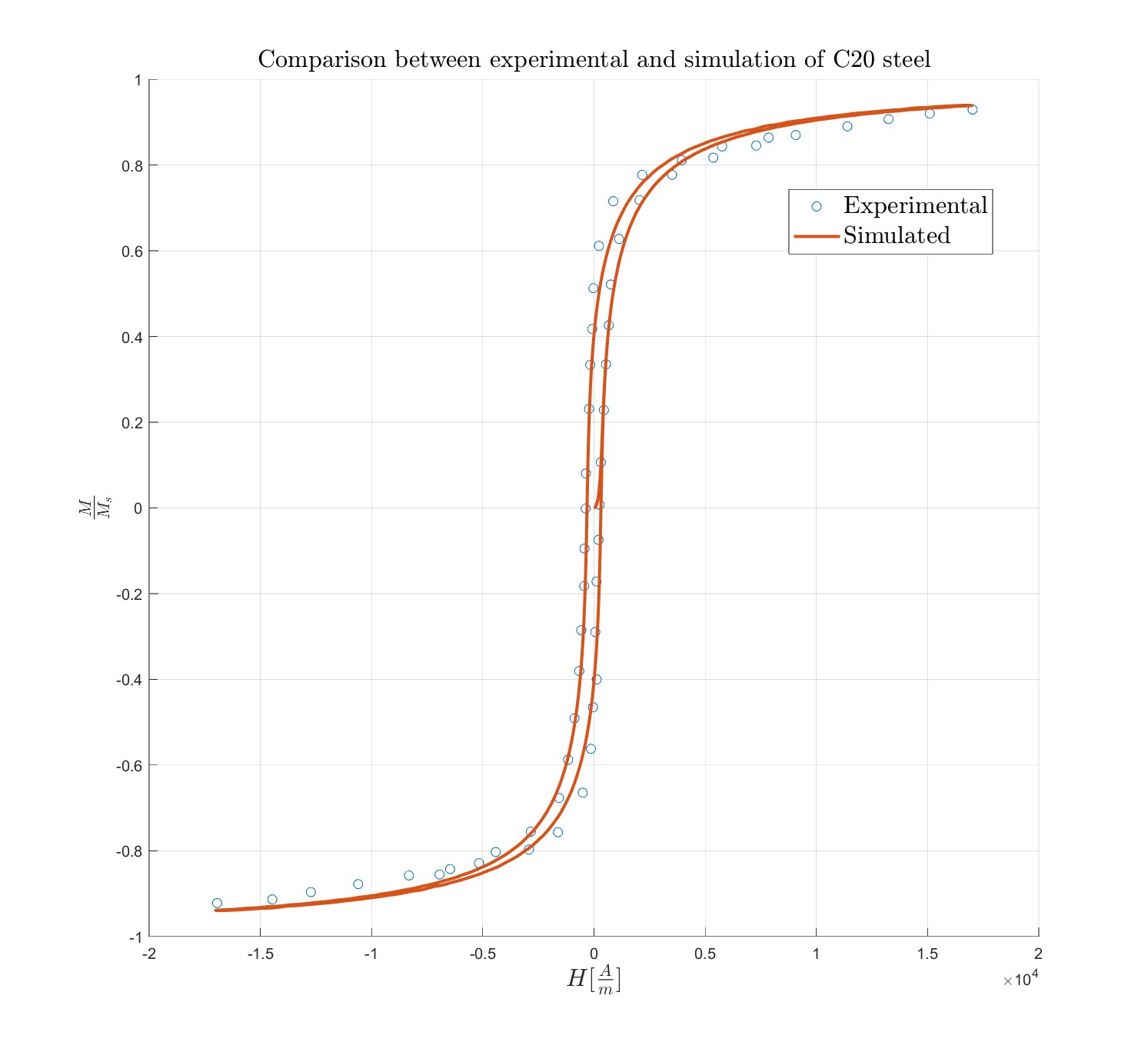}
\caption{\label{fig:IMG10_11} Experimental hysteresis loop and its simulation from \cite{Jiles1992}.}
\end{figure}


In the case of real data, a hysteresis loop is obtained from a Non Oriented M 470-50A produced by Marcegaglia Ravenna s.p.a with a Single Sheet Tester from Brockhaus Messtechnik. This machine has the following characteristics:
\begin{itemize}
\item Model: MPG100 D DC/AC
\item frequency ranges: from 3Hz to 10 kHz
\item maximum polarization: 2T
\item measurement repeatability: $\leq 2 $ percent;
\end{itemize}

From the result represented in figure \ref{fig:IMG12}, we can see is a good agreement between experimental data points and simulation.
\begin{figure}[h!]
\centering
\includegraphics[width=0.8\textwidth]{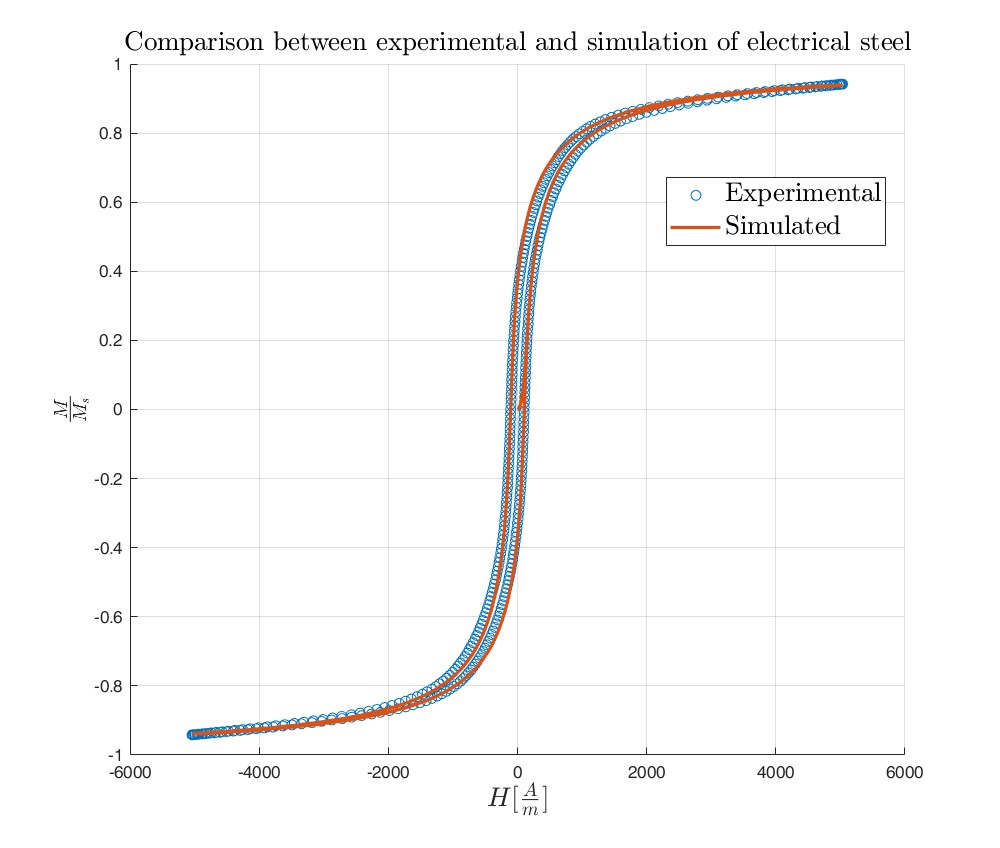}
\caption{\label{fig:IMG12} Experimental hysteresis loop and its simulation.}
\end{figure}
\section{Conclusion}
This paper focused on the Jiles-Atherton Model, which is widely used in engineering applications, and presented a new approach for finding the simulation parameters for the anhysteretic curve of ferromagnetic materials. By using the material's susceptibility and linearizing the anhysteretic magnetization curve with a paramagnetic function, we could find the magnetic moments of the material and determine the simulation parameters in a more physical and simplified manner. Our results showed that it is possible to describe the anhysteretic magnetization curve of a ferromagnetic material with a linear approximation of a paramagnet for every value of the external applied field. 
Validation of the proposed method with synthetic and experimental data has demonstrated its effectiveness and stability. \\
In conclusion, \name{} extends the approach of Algorithm \ref{alg:estimation} by improving the quality of parameter estimation without requiring the iterative solution of a system of ordinary differential equations (ODEs), which is computationally expensive and presents challenges in solving an inverse problem.
This approach can be useful in many engineering applications requiring accurate ferromagnetic material characterisation.
\section{Acknowledgement }
This work was financed by the European Union - NextGenerationEU (National Sustainable Mobility Center CN00000023, Italian Ministry of University and Research Decree n. 1033 - 17/06/2022, Spoke 11 - Innovative Materials \& Lightweighting), and National Recovery and Resilience Plan (NRRP), Mission 04 Component 2 Investment 1.5 – NextGenerationEU, Call for tender n. 3277 dated 30/12/2021.
Moreover this work was supported by Alessandro Ferraiuolo, R \& D Manager of Marcegaglia Ravenna S.p.A., giving the material for experimental validations.

\end{flushleft}

\end{document}